\documentclass[a4paper,12pt]{article}
\usepackage{amsmath}
\usepackage{hyperref}
\usepackage{graphicx}

\input{tcilatex}
\begin{document}

\title{Determining the Rolle function in Lagrange interpolatory approximation%
}
\author{J. S. C. Prentice \\
Senior Research Officer\\
Mathsophical\\
Johannesburg, South Africa}
\maketitle

\begin{abstract}
We determine the Rolle function in Lagrange polynomial approximation using a
suitable differential equation. We then propose a device for improving the
Lagrange approximation by exploiting our knowledge of the Rolle function.
\end{abstract}

\section{Introduction}

Approximation of nonlinear functions is of fundamental importance in applied
mathematics, especially by means of polynomials. Those techniques for which
the approximation error is well understood are particularly useful. In this
paper, we describe how the error in Lagrange polynomial interpolation can be
precisely determined, by solving an appropriate initial-value problem. We
then consider using the knowledge so obtained to improve the quality of the
original approximation.

\section{Relevant Concepts, Terminology and Notation}

Let $f\left( x\right) $ be a real-valued univariate function. The \textit{%
Lagrange interpolating polynomial} $P_{n}\left( x\right) $ of degree $n,$ at
most, that interpolates the data $\{f\left( x_{0}\right) ,$ $f\left(
x_{1}\right) ,\ldots ,f\left( x_{n}\right) \}$ at the nodes $\left\{
x_{0},x_{1},\ldots ,x_{n}\right\} ,$ where $x_{0}<x_{1}<\cdots <x_{n},$ has
the property%
\begin{equation*}
P_{n}\left( x_{k}\right) =f\left( x_{k}\right)
\end{equation*}%
for $k=0,1,\ldots ,n.$ Naturally, we regard $P_{n}\left( x\right) $ as an
approximation to $f\left( x\right) .$ The pointwise error in Lagrange
interpolation, on $\left[ x_{0},x_{n}\right] ,$ is 
\begin{equation}
\Delta \left( x;P_{n}\right) \equiv f\left( x\right) -P_{n}\left( x\right) =%
\frac{f^{\left( n+1\right) }\left( \xi \left( x\right) \right) }{\left(
n+1\right) !}\dprod\limits_{k=0}^{n}\left( x-x_{k}\right) ,  \label{delta}
\end{equation}%
where $x_{0}<\xi \left( x\right) <x_{n},$ and is derived by invoking Rolle's
Theorem \cite{I and K}\cite{K and C}. Clearly, we assume here that $f\left(
x\right) $ is $\left( n+1\right) $-times differentiable and, as will be
seen, we must assume that $f\left( x\right) $ is, in fact, $\left(
n+2\right) $-times differentiable. We refer to $\xi $ generically as the 
\textit{Rolle number}, and to $\xi \left( x\right) $ as the \textit{Rolle
function}.

\section{Determining the Rolle Function}

Using the notation $Q_{n}\left( x\right) \equiv
\dprod\nolimits_{k=0}^{n}\left( x-x_{k}\right) $ for the remainder
polynomial, we have, by differentiating with respect to $x,$%
\begin{align*}
\left( n+1\right) !\left( f\left( x\right) -P_{n}\left( x\right) \right) &
=f^{\left( n+1\right) }\left( \xi \left( x\right) \right) Q_{n}\left(
x\right) \\
\Rightarrow \left( n+1\right) !\left( f^{\prime }\left( x\right)
-P_{n}^{\prime }\left( x\right) \right) & =f^{\left( n+1\right) }\left( \xi
\right) Q_{n}^{\prime }\left( x\right) +Q_{n}\left( x\right) \frac{%
df^{\left( n+1\right) }\left( \xi \right) }{d\xi }\frac{d\xi }{dx} \\
& =f^{\left( n+1\right) }\left( \xi \right) Q_{n}^{\prime }\left( x\right)
+Q_{n}\left( x\right) f^{\left( n+2\right) }\left( \xi \right) \frac{d\xi }{%
dx}.
\end{align*}%
We have used the well-known prime notation for differentiation with respect
to $x.$ In this expression, the factor $f^{\left( n+1\right) }\left( \xi
\right) $ denotes the $\left( n+1\right) $th derivative of $f\left( \xi
\right) $ with respect to $\xi ,$ and similarly for $f^{\left( n+2\right)
}\left( \xi \right) .$ We now find%
\begin{equation*}
\frac{d\xi }{dx}=\frac{\left( n+1\right) !\left( f^{\prime }\left( x\right)
-P_{n}^{\prime }\left( x\right) \right) -f^{\left( n+1\right) }\left( \xi
\right) Q_{n}^{\prime }\left( x\right) }{Q_{n}\left( x\right) f^{\left(
n+2\right) }\left( \xi \right) }
\end{equation*}%
and, if we have a particular value $\xi _{z}=\xi \left( x_{z}\right) $ at
our disposal, we then have an initial-value problem that can, in principle,
be solved to yield the Rolle function $\xi \left( x\right) .$ Note the
necessity of our assumption that $f\left( x\right) $ is $\left( n+2\right) $%
-times differentiable.

\section{Calculations}

Consider the Lagrange interpolation of%
\begin{equation*}
f\left( x\right) =e^{x}\sin x
\end{equation*}%
over the nodes $\left\{ 0,\frac{3\pi }{2}\right\} .$ We have $n=1$ so that
the Lagrange polynomial is 
\begin{equation*}
P_{1}\left( x\right) =\left( \frac{e^{\frac{3\pi }{2}}\sin \left( \frac{3\pi 
}{2}\right) }{\frac{3\pi }{2}}\right) x=\left( -\frac{2e^{\frac{3\pi }{2}}}{%
3\pi }\right) x.
\end{equation*}%
Furthermore, we have%
\begin{equation}
\Delta \left( x;P_{1}\right) =e^{x}\sin x-\left( -\frac{2e^{\frac{3\pi }{2}}%
}{3\pi }\right) x=\left( x^{2}-\frac{3\pi x}{2}\right) e^{\xi \left(
x\right) }\cos \left( \xi \left( x\right) \right)  \label{delta example 1}
\end{equation}%
and%
\begin{equation}
\frac{d\xi }{dx}=\frac{2\left( e^{x}\left( \cos x+\sin x\right) +\frac{2e^{%
\frac{3\pi }{2}}}{3\pi }\right) -\left( 2e^{\xi }\cos \xi \right) \left( 2x-%
\frac{3\pi }{2}\right) }{\left( x^{2}-\frac{3\pi x}{2}\right) \left( 2e^{\xi
}\left( \cos \xi -\sin \xi \right) \right) }.  \label{DE example 1}
\end{equation}%
We solve this differential equation using an initial value chosen close to
the node $x_{0}=0.$ Observe that it is not possible to find the Rolle number
at any interpolation node, because the factor $\dprod\nolimits_{k=0}^{n}%
\left( x-x_{k}\right) $ in (\ref{delta}) ensures that $\Delta \left(
x;P_{n}\right) =0$ at the interpolation nodes, \textit{irrespective of the
value of} $\xi .$ Hence, we choose here $x_{z}=10^{-5}$ and determine the
corresponding Rolle number by applying Newton's method to (\ref{delta
example 1}). Naturally, we compute $\Delta \left( x_{z};P_{1}\right)
=e^{x_{z}}\sin x_{z}+\left( \frac{2e^{\frac{3\pi }{2}}}{3\pi }\right) x_{z}$
to facilitate this calculation. In fact, we find two values for $\xi _{z},$
giving the initial values $\left( x_{z},\xi _{z}\right) =\left(
10^{-5},2.1931\right) $ and $\left( x_{z},\xi _{z}\right) =\left(
10^{-5},4.6631\right) .$ (Note to the reader: we quote numerical values
correct to no more than four decimal places throughout this paper, simply
for ease of presentation, but all calculations are performed in double
precision). Using these two initial values we solve (\ref{DE example 1}) to
find $\xi \left( x\right) ,$ shown respectively in Figures 1(a) and 1(b). We
are then able to compute the pointwise error $\Delta \left( x;P_{1}\right) $
using (\ref{delta example 1}), and this is shown in Figure 1(c) for $\left(
x_{z},\xi _{z}\right) =\left( 10^{-5},2.1931\right) .$ Of course, we can
compare this error curve with the actual error, and the magnitude of the
difference between the two - the `error in $\Delta $', so to speak - is
shown in Figure 1(d). Clearly, this error is small, indicating the accuracy
of our numerical estimate of $\xi \left( x\right) .$ Similar results obtain
for $\left( x_{z},\xi _{z}\right) =\left( 10^{-5},4.6631\right) ,$ with a
maximum magnitude in the error in $\Delta $ of $\sim 2\times 10^{-11}.$ This
accuracy is a consequence of using a fourth-order explicit Runge-Kutta
method \cite{Butcher}\cite{Hairer et al} to solve (\ref{DE example 1}) with
a small stepsize $\left( \sim 5\times 10^{-5}\right) .$

As a second example, let us use the same objective function $f\left(
x\right) $ as above, but now with three interpolatory nodes $\left\{ 0,2,%
\frac{3\pi }{2}\right\} .$ So, we have $n=2$ which gives 
\begin{equation*}
P_{2}\left( x\right) =ax^{2}+bx+c,
\end{equation*}%
where $a=-9.9476,b=23.2546$ and $c=0.$ Also, with 
\begin{align*}
Q_{2}\left( x\right) & =\left( x-0\right) \left( x-2\right) \left( x-\frac{%
3\pi }{2}\right) \\
& =\left( x^{3}+\left( -\frac{3\pi }{2}-2\right) x^{2}+3\pi x\right) \\
Q_{2}^{\prime }\left( x\right) & =3x^{2}+\left( -3\pi -4\right) x+3\pi ,
\end{align*}%
we have%
\begin{equation*}
\Delta \left( x;P_{2}\right) =\frac{e^{\xi }\left( \cos \xi -\sin \xi
\right) Q_{2}\left( x\right) }{3}
\end{equation*}%
and%
\begin{equation}
\frac{d\xi }{dx}=\frac{6\left( e^{x}\left( \cos x+\sin x\right) -\left(
2ax+b\right) \right) -2e^{\xi }\left( \cos \xi -\sin \xi \right)
Q_{2}^{\prime }\left( x\right) }{Q_{2}\left( x\right) \left( -4e^{\xi }\sin
\xi \right) }.  \label{DE example 2}
\end{equation}%
We perform similar calculations as above, with results depicted in Figure 2.
Again, we find two initial values $\left( x_{z},\xi _{z}\right) =\left(
10^{-5},1.7845\right) $ and $\left( x_{z},\xi _{z}\right) =\left(
10^{-5},3.8165\right) ,$ and the error plot in Figure 2(d) corresponds to $%
\left( x_{z},\xi _{z}\right) =\left( 10^{-5},1.7845\right) .$ For $\left(
x_{z},\xi _{z}\right) =\left( 10^{-5},3.8165\right) $ the maximum magnitude
in the error in $\Delta $ is $\sim 3\times 10^{-12}.$

Our second example allows us to make a point about the continuity of $\xi
\left( x\right) .$ Since $P_{2}\left( x\right) $ is continuous and the
assumed $\left( n+2\right) $-times differentiability of $f\left( x\right) $
implies the continuity of both $f\left( x\right) $ and $f^{\left( n+1\right)
}\left( x\right) ,$ we will assume, from (\ref{delta}), that $\xi \left(
x\right) $ is continuous. Hence, even though it is not really meaningful to
ask for the value of $\xi \left( x\right) $ at an interpolatory node, the
continuity of $\xi \left( x\right) $ allows us to infer a value at such a
node. We refer to such a value as an \textit{implied} Rolle number. In our
second example we determine the implied values $\xi \left( 2\right) =2.0991$
from Figure 2(a) and $\xi \left( 2\right) =3.7381$ from Figure 2(b). Our
technique for doing this is as follows: we ensure that $x=2$ is not amongst
the nodes used for the Runge-Kutta calculation, because this would lead to a
zero in the denominator on the RHS of (\ref{DE example 2}), but we use the
values of $\xi \left( x\right) $ at the Runge-Kutta nodes on either side of $%
x=2$ to estimate $\xi \left( 2\right) $ using linear interpolation. Of
course, we acknowledge that it is not necessary to know the Rolle number at
any interpolatory node, because the pointwise error at interpolatory nodes
is always zero. As such, $\xi \left( x\right) $ is actually arbitrary at
interpolatory nodes, but our technique does allow a sensible estimate to be
made, if only for completeness' sake. Note also that simple linear
extrapolation will allow an estimate of $\xi \left( x\right) $ to be made at
the endpoints of the interval.

\section{A\ Possible Application}

Our ability to determine the Rolle function suggests an interesting
possibility. If we know $\xi \left( x\right) $ then we know $f^{\left(
n+1\right) }\left( \xi \left( x\right) \right) .$ If we approximate $%
f^{\left( n+1\right) }\left( \xi \left( x\right) \right) $ by means of a
polynomial - a least-squares fit, or a truncated Taylor series, for example
- then, using (\ref{delta}), we find%
\begin{equation*}
f\left( x\right) \approx P_{n}\left( x\right) +\frac{P_{\xi }\left( x\right) 
}{\left( n+1\right) !}\dprod\limits_{k=0}^{n}\left( x-x_{k}\right) ,
\end{equation*}%
where $P_{\xi }\left( x\right) $ denotes the polynomial that approximates $%
f^{\left( n+1\right) }\left( \xi \left( x\right) \right) .$ Notice that the
RHS of this expression is simply a polynomial, and so constitutes a
polynomial approximation to $f\left( x\right) .$ Thus, our knowledge of $\xi
\left( x\right) $ allows us to improve the approximation $P_{n}\left(
x\right) $ by adding a polynomial term that approximates the error in $%
P_{n}\left( x\right) .$

By way of example, we return to the first of our earlier calculations. Here,
we have 
\begin{eqnarray*}
f\left( x\right) &=&e^{x}\sin x \\
P_{1}\left( x\right) &=&\left( -\frac{2e^{\frac{3\pi }{2}}}{3\pi }\right) x
\\
\Delta \left( x;P_{1}\right) &=&\left( x^{2}-\frac{3\pi x}{2}\right) e^{\xi
\left( x\right) }\cos \left( \xi \left( x\right) \right)
\end{eqnarray*}%
using the nodes $\left\{ 0,\frac{3\pi }{2}\right\} .$ Once we have found $%
\xi \left( x\right) $ - using $\left( x_{z},\xi _{z}\right) =(10^{-5},$ $%
2.1931)$ and shown in Figure 1(a) - we determine a polynomial approximation
to $e^{\xi \left( x\right) }\cos \left( \xi \left( x\right) \right) $ using
a least-squares fit. For the purpose of demonstration, we choose to fit a
degree six polynomial so that $\Delta \left( x;P_{1}\right) +P_{1}\left(
x\right) $ becomes an eighth-degree polynomial approximation, which we
denote $P_{8}^{\Delta }\left( x\right) .$ In Figure 3(a) we show the
pointwise error in $P_{1}\left( x\right) ,$ and observe it has a maximum
magnitude of $\sim 75.$ In Figure 3(b) we plot the pointwise error 
\begin{equation*}
e^{x}\sin x-P_{8}^{\Delta }\left( x\right) .
\end{equation*}%
This error is significantly smaller than that in Figure 3(a) and has maximum
magnitude $\sim 6\times 10^{-3}.$ In other words, $P_{8}^{\Delta }\left(
x\right) $ is an approximation more than four orders of magnitude (!) more
accurate than $P_{1}\left( x\right) ,$ and this was achieved only through
our knowledge of the Rolle function $\xi \left( x\right) .$ This example
serves to illustrate the value of knowing $\xi \left( x\right) ,$ and
certainly warrants further investigation, but we will reserve a detailed
study thereof for future research.

\section{Conclusion}

We have described how the Rolle function in Lagrange interpolatory
polynomial approximation can be determined by solving an initial-value
problem. Knowledge of the Rolle function permits the calculation of the
approximation error. In particular, the Rolle term in the expression for the
approximation error can itself be approximated by means of a polynomial,
once the Rolle function is known, and this can lead to a significant
improvement in the quality of the Lagrange approximation overall. Of course,
the ideas presented in this paper are not restricted to Lagrange
approximation, but could also be applied to Hermite interpolation, for
example, in which the error term is also derived through the use of Rolle's
Theorem.

\begin{equation*}
\end{equation*}

\begin{figure}[tbp]
\begin{center}
\includegraphics[width=1\linewidth]{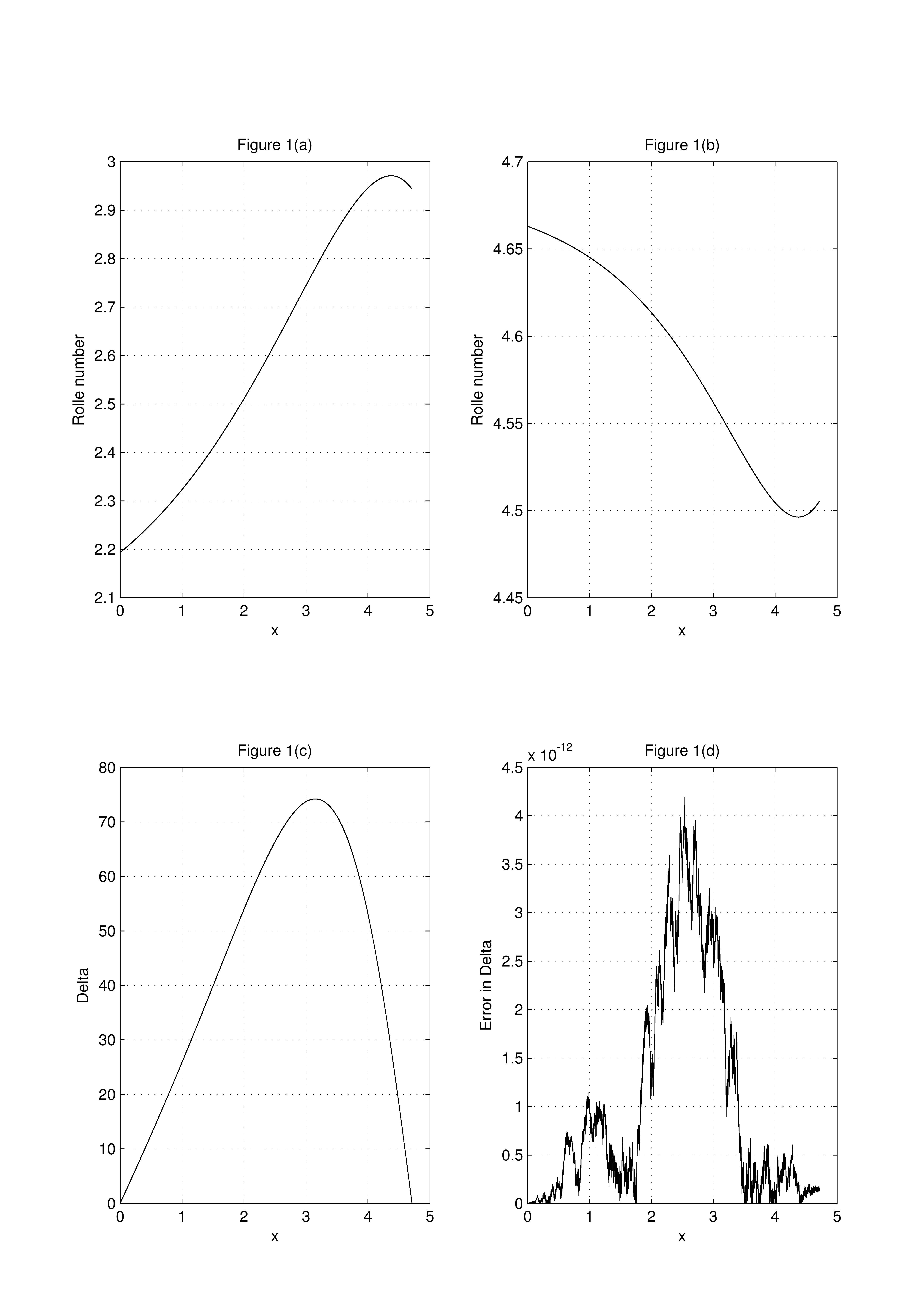} 
\end{center}
\end{figure}
\newpage

\begin{figure}[tbp]
\begin{center}
\includegraphics[width=1\linewidth]{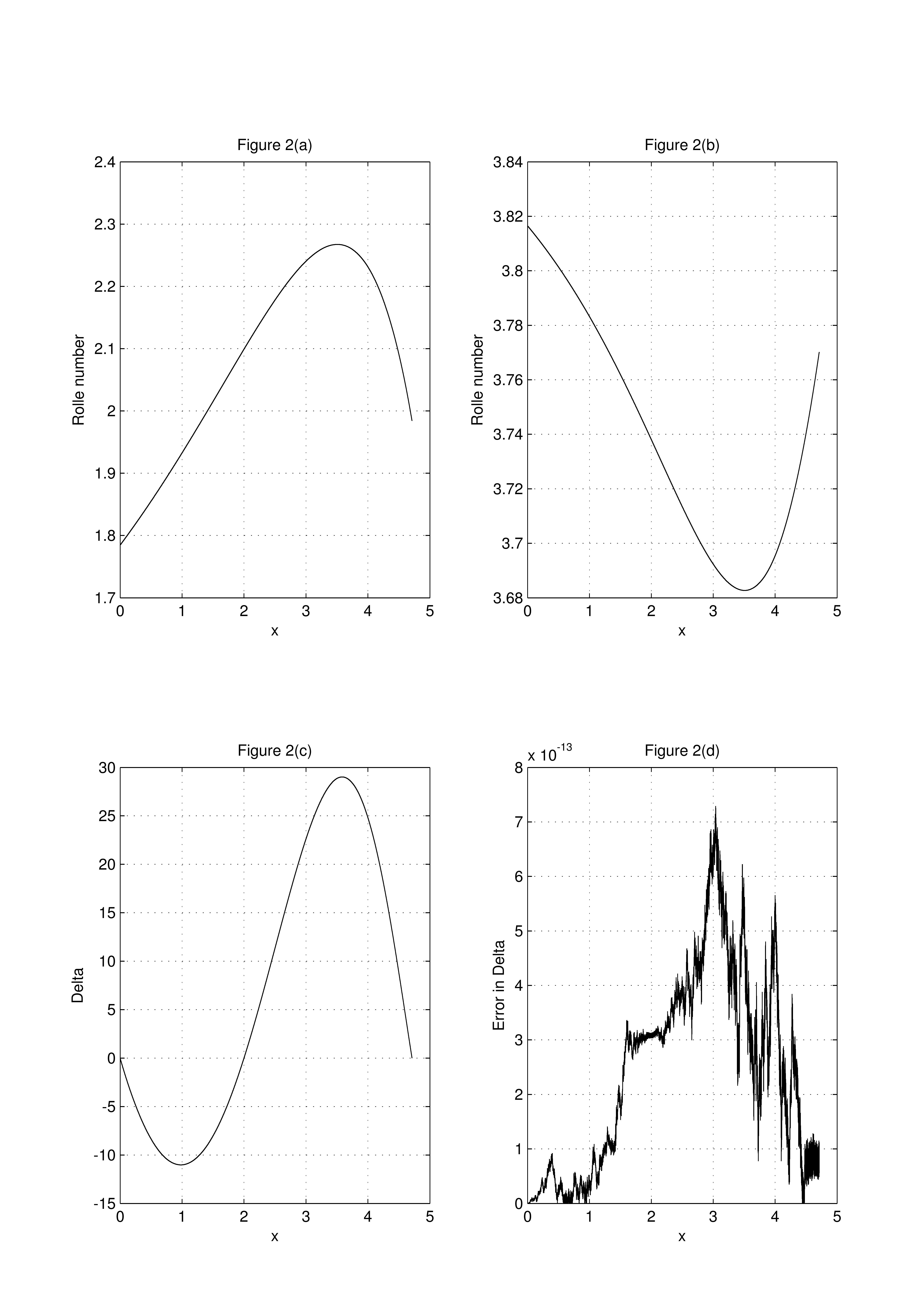} 
\end{center}
\end{figure}
\newpage

\begin{figure}[tbp]
\begin{center}
\includegraphics[width=1\linewidth]{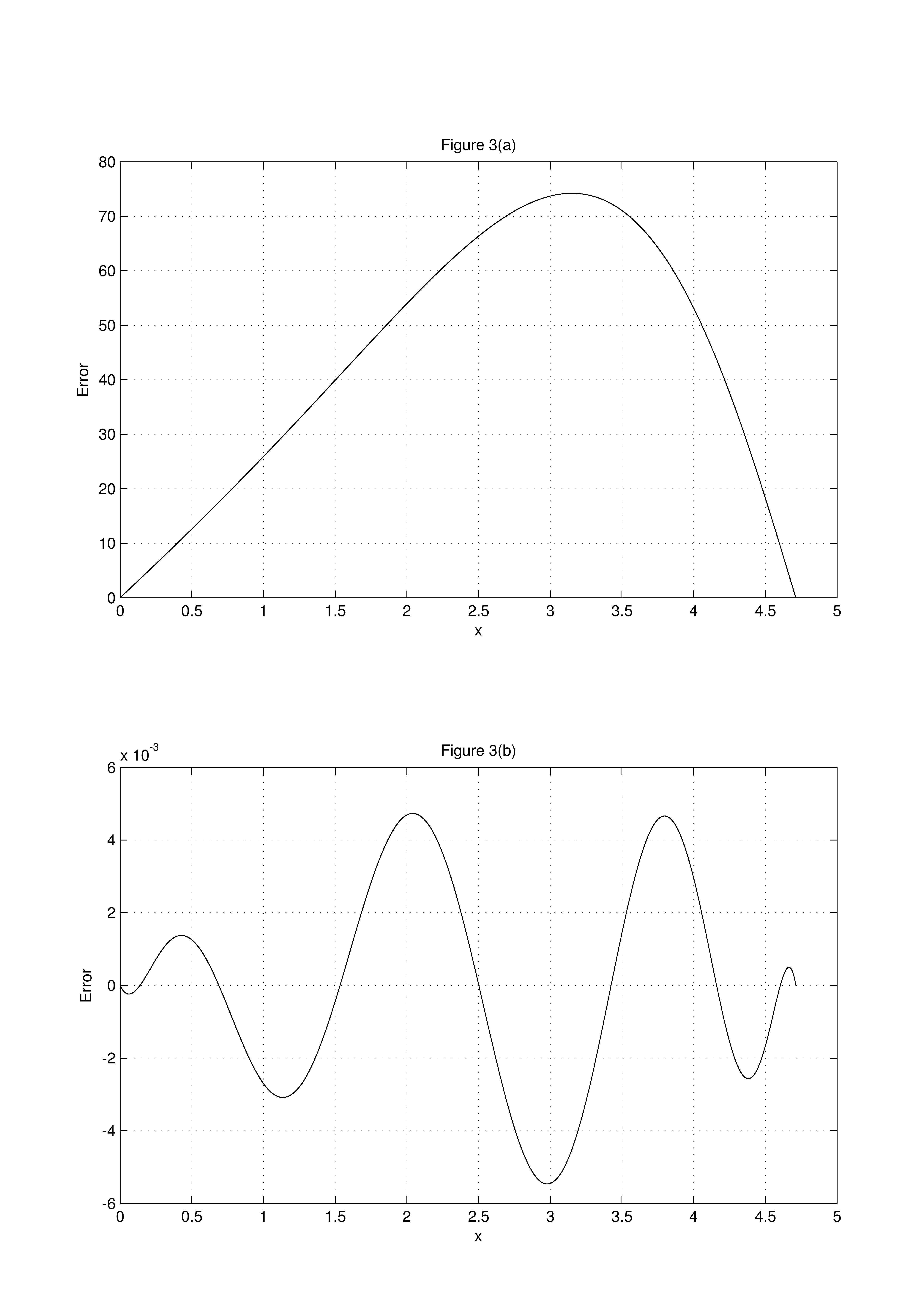} 
\end{center}
\end{figure}

\end{document}